\documentclass{article}
\usepackage{amssymb}
\newtheorem{theo}{\arabic{section}.\arabic{abz}. Theorem}

\newtheorem{defi}{\arabic{section}.\arabic{abz}. Definition}

\newtheorem{propo}{\arabic{section}.\arabic{abz}. Proposition}

\newtheorem{coro}{\arabic{section}.\arabic{abz}. Corollary}

\newtheorem{lemm}{\arabic{section}.\arabic{abz}. Lemma}

\newtheorem{rema}{\arabic{section}.\arabic{abz}. Remark}

\newtheorem{exa}{\arabic{section}.\arabic{abz}. Example}

\newcounter{abz}[section]
\newcounter{sabz}[abz]
\addtocounter{section}{-1}
\newcommand{\abz}{\refstepcounter{abz}}
\newcommand{\sabz}{\refstepcounter{sabz}}

\def\Ad{\mathop{\mathrm{Ad}}}

\def\Aut{\mathop{\mathrm{Aut}}}
\def\codim{\mathop{\mathrm{codim}}}
\def\crk{\mathop{\mathrm{corank}}}
\def\diag{\mathop{\mathrm{diag}}}
\def\End{\mathop{\mathrm{End}}}
\def\gl{\mathop{\mathfrak{gl}}}

\def\Id{\mathop{\mathrm{Id}}}
\def\im{\mathop{\mathrm{im}}}

\def\ind{\mathop{\mathrm{ind}}}

\def\rk{\mathop{\mathrm{rank}}}
\def\Sing{\mathop{\mathrm{Sing}}}
\def\sl{\mathop{\mathfrak{sl}}}
\def\so{\mathop{\mathfrak{so}}}

\def\Sp{\mathop{\mathrm{Sp}}}
\def\Tr{\mathop{\mathrm{Tr}}}

\def\b{\mathfrak{b}}
\def\C{\mathbb{C}}
\def\const{\mathrm{const}}
\def\D{{\mathfrak g}^*}

\def\F{\mathbb{K}}

\def\g{{\mathfrak g}}

\def\h{{\mathfrak h}}

\def\k{{\mathfrak k}}

\def\n{\mathfrak{n}}

\def\R{\mathbb{R}}

\def\t{\theta}

\def\la{\lambda}

\def\sdp{\rhd\!\!\!<}
\def\qed{Q.E.D.}
\def\pf{\noindent{\sc Proof}\ }
\renewcommand{\le}{\leqslant}
\renewcommand{\ge}{\geqslant}
\renewcommand{\phi}{\varphi}

\begin{document}
\title{Algebraic Nijenhuis operators and Kronecker Poisson pencils}
\author{Andriy
Panasyuk}

\maketitle
\bibliographystyle{amsalpha}
\section{Introduction}
This paper is devoted to a method of constructing completely integrable systems based on the micro-local theory of bihamiltonian structures \cite{gz1,gz2,bols,gz3,gz4,p1,z1}.
The main tool are the so-called micro-Kronecker bihamiltonian structures \cite{z1}, which will be called Kronecker in this paper for short (in \cite{gz4} the term Kronecker was used for the micro-Kronecker structures with some additional condition of "flatness" which will not be essential in this paper).

A Kronecker bihamiltonian structure on a manifold $M$ is a Poisson pencil $\{s_1\t_1+s_2\t_2\}_{(s_1,s_2)\in\F^2}$, i. e. a two-dimensional linear space over a base field $\F$ in the set of all Poisson structures on $M$, satisfying an additional condition of the constancy of rank: $\rk_\C\t^s=\const,s:=(s_1,s_2)\in\C^2\setminus(0,0),
\t^s:=s_1\t_1+s_2\t_2$    (in the real case we should pass to the complexification of the pencil). The kroneckerity condition is important due to the fact that it automatically implies the existence (at least locally) on $M$ of the complete involutive with respect to any bivector $\t^s$
set of functions. This set is functionally generated by the Casimir functions  of the  bivectors $\t^s$ (see Proposition \ref{20.30}). Geometrically this set corresponds to the intersection over $s\in\F^2\setminus\{(0,0)\}$ of all symplectic leaves of maximal dimension of the Poisson structures $\t^s$ and the completeness of this set reflects the fact that this intersection is lagrangian in any fixed symplectic leaf (see \cite{gz2,p1}).

The main result of this paper (Theorem \ref{main1}) gives a criterion of kroneckerity for the Poisson pencils related to  diagonalizable algebraic Nijenhuis operators. An algebraic Nijenhuis operator $N$ on a Lie algebra $\g$ (see \cite{mks,cgm1}, for example) is a linear operator $N:\g\to \g$ with the condition of the vanishing of the so-called Nijenhuis torsion (see Definition \ref{10.10}). Given a linear operator $N:(\g,[,])\to (\g,[,])$, the condition of vanishing of its Nijenhuis torsion guarantees that the infinithesimal part $[,]_N$ of the trivial deformation $(\Id+\la N)^{-1}[(\Id+\la N)\cdot,(\Id+\la N)\cdot]$ of the Lie bracket $[,]$ is again a Lie bracket. This new Lie bracket $[,]_N$ is automatically compatible with $[,]$, thus any Nijenhuis operator $N$ "produces" the pencil of Lie brackets $[,]^s:=s_1[,]+s_2[,]_N$ and, consequently, the corresponding pencil $\{\t_N^s\}_{s\in\F^2}$ of the Lie-Poisson structures on $\g^*$. 

Let us look more closely at the problem of the kroneckerity of the Poisson pencil $\{\t_N^s\}$. It can be shown (Proposition \ref{10.20}) that if $N$ is Nijenhuis, $(N-\la\Id)^{-1}[(N-\la\Id)\cdot,(N-\la\Id)\cdot]=[\cdot,\cdot]_N
-\la[\cdot,\cdot]$. In particular, all the Lie brackets $[,]^s$ are isomorphic to $[,]$ except those corresponding to $s=(s_1,s_2)$ with $\la=-s_1/s_2$ belonging to $\Sp N$, the spectrum of $N$. Thus the problem of kroneckerity of $\{\t_N^s\}$ (modulo some not very restrictive assumption on the codimension of the set of singular coadjoint orbits of $\g$, see (\ref{e.main})) reduces to the problem of calculating the dimension of the coadjoint orbits of the exceptional brackets $[,]_N-\la_i[,], i=1,\ldots,n$, where $\la_1,\ldots,\la_n$ are the eigenvalues of $N$. In fact, due to the semicontinuity of the function $\rk\t^s$, in order to prove the kroneckerity it is sufficient to find for any $i$ a particular coadjoint orbit $O_i$ of a Lie bracket $[,]_N-\la_i[,]$ such that $\dim O_i=\dim O$, where $O$ is the generic coadjoint orbit of $[,]$.
 One possibility of finding the orbits $O_i$ is the following. If $N$ is a Nijenhuis operator, then $N: (\g,[,]_N)\to (\g,[,])$ is a homomorphism of Lie algebras \cite{mks}. Hence $\im N$ is a subalgebra of $(\g,[,])$ and we have a Poisson inclusion $^t\! N:((\im N)^*,\t_{\mathrm{st}})\hookrightarrow(\g^*,\t_N)$, where $\t_{\mathrm{st}}$ is the standard Lie--Poisson structure on $(\im N)^*$ and $\t_N$ corresponds to $[,]_N$. In particular, one can take $O_i$ to be a symplectic leaf in $([\im (N-\la_i\Id)]^*,\t_{\mathrm{st}})\subset(\g^*,\t_{(N-\la_i\Id)})$ (the operator $N-\la_i\Id$ is also  Nijenhuis). Choosing $O_i$ to be a generic coadjoint orbit and passing to codimensions we get the following sufficient condition of kroneckerity: if $\ind\im (N-\la_i\Id)+\codim\im (N-\la_i\Id)=\ind\g$ for any $i$, where $\ind$ stands for the index of a Lie algebra, i. e. the codimension of a generic coadjoint orbit, then the Poisson pencil $\{\t_N^s\}$ is Kronecker (cf. Corollary \ref{coro}).
 
 In general, however, this condition is not necessary because it may happen that the generic coadjoint orbits in $(\im N)^*$ are not generic in $(\g^*,\t_N)$. For example, take $\g=\sl(2)=\n_-\oplus\b_+$, where $\b_+$ is the upper Borel subalgebra and $\n_-$ is the lower nilpotent subalgebra. Let $N=P_{\n_-}$ be the projector to the first summand along the second one. Then coadjoint orbits of $\im N$ are points, whereas the algebra $(\g,[,]_N)$ is nonabelian and has also coadjoint orbits of dimension 2. 

So our main theorem generalizes the above mentioned sufficient condition and gives necessary and sufficient conditions of the kroneckerity of the pencil $\{\t^s_N\}$ (for a diagonalizable $N$). The method of proof of this result consists in showing that the above mentioned exceptional brackets are in fact semi-direct products and using the Ra\"{\i}s type formulas for their indices.

We illustrate our method by two examples. First of them, a generalization of the example above, relates a complete involutive set of functions on a semisimple split Lie algebra with the Nijenhuis operator $N$ that is a projector onto the lower nilpotent subalgebra along the upper Borel subalgebra (in fact we use operators of the form $s_1N+s_2\Id$, all such operators generate the same Poisson pencil,  see Section \ref{s.borel}).
 
Note that such $N$ is a classical $R$-matrix in the sense of Semenov-Tian-Shansky \cite{sts}. However, our method is essentially different, since it 1) uses another modified bracket $[\cdot,\cdot]_N=[N \cdot,\cdot]+[\cdot,N \cdot]-N[\cdot,\cdot]$ (in the $R$-matrix approach $[\cdot,\cdot]_R=[R \cdot,\cdot]+[\cdot,R \cdot]$); 2) uses the whole pencil of brackets, generated by $[\cdot,\cdot]$ and $[\cdot,\cdot]_N$ (the  $R$-matrix approach uses only $[\cdot,\cdot]_R$); 3) is applicable to the generic coadjoint orbits (while $R$-matrix approach generates involutive sets of functions on the orbits of dimension $2r=2\ind \g$). The reader is also referred to the reference \cite{mks} for another application of algebraic Nijenhuis operators which are projectors.  This application, related with the so-called  Kostant--Symes theorem, is also different from our since it produces involutive sets of functions on the dual spaces to subalgebras of $\g$.

Our second example (see Section \ref{s.manakov}), which in fact inspired this paper, uses the Lie algebra $\g=\gl_n$ of $n\times n$-matrices and the operator $N=L_A$ of the left multiplication by a diagonal matrix $A=\diag(\la_1,
\ldots,\la_n)\in\g$
with $\la_i\not=\la_j,i\not=j$. The corresponding modified commutator is of the form $[x,y]_N=xAy-yAx, x,y\in\g$, and was considered earlier in \cite{bols} and \cite{morosip} in the context of the  free rigid body system on $\so_n$. In these papers it was proved that the related Euler vector field is hamiltonian with respect to the corresponding  Lie--Poisson structure $\t_N$. It is  also known that this vector field is hamiltonian with respect to the standard Lie--Poisson structure $\t_\mathrm{st}$ on $\so_n$ with the hamiltonian function $\Tr(Ax^2)$. As a consequence of our method we get an alternative proof of  the complete integrability of the analogue of the $n$-dimensional free rigid body system on $\gl_n$. The traditional proof, which goes back to the papers of Manakov \cite{man} and Mishchenko--Fomenko \cite{mf} (see also \cite{bols}), uses the so-called method of the argument translation, i.e. the pencil of the affine Poisson structures generated by the linear Poisson structure $\t_\mathrm{st}$ and the constant Poisson structure $\t_\mathrm{st}|_A$. The complete involutive family of functions (which includes the function $\Tr(Ax^2)$) built by our method in fact coincides with the family obtained by the method of the argument translation (see Proposition \ref{manakov}), however the two families of functions are obtained differently and their equality is not seen at first glance. The proof of the equality uses  recurrence relations between two families. 

In order to use our method for the proof of the integrability of "true" free rigid body (i. e. on $\so_n$) one should extend the method. Note that the pencil of brackets $[x,y]_{(N-\la\Id)}=x(A-\la\Id)y-y(A-\la\Id)x$ is correctly defined on the subalgebra $\so_n \subset\gl_n$, although the Nijenhuis operator $N$ does not preserve this sublagebra: $N\so_n\not \subset\so_n$ (cf. Example \ref{10.60} and Remark \ref{10.70}).
One can consider the following generalization of this situation: let $\g$ be a Lie algebra, $\k \subset \g$ its subalgebra, $N:\g\to\g$ a Nijenhuis operator such that $(N-\la\Id)\k$ again is a Lie subalgebra for any $\la$ (we don't require $N\k \subset\k$). Then $(N-\la\Id)^{-1}[(N-\la\Id)\cdot,(N-\la\Id)\cdot]$ is a Lie bracket on $\k$ which is equal to $[\cdot,\cdot]_N|_\k-\la[\cdot,\cdot]|_\k=:[\cdot,\cdot]^\la$.

 So one of the possible extensions of our method is the study of the kroneckerity of the  pencils of the Lie--Poisson structures related to the pencil of the brackets $[\cdot,\cdot]^\la$. Another one is the consideration of the so-called weak Nijenhuis operators \cite{cgm1}, i. e. operators whose Nijenhuis torsion is a cocycle with the coefficients in the adjoint module. Such operators also generate pencils of Lie--Poisson structures and the question of their kroneckerity  seems reasonable and can provide with new examples of completely integrable systems or new proofs of their complete integrability.

\section{Algebraic Nijenhuis operators and pencils of Lie algebras}
\label{s10} The following definition, which is basic for this
paper, is taken from \cite{mks}. 
\abz\label{10.10}
\begin{defi}
\rm Let $(\g,[,])$ be a finite-dimensional Lie algebra over a
field $\F$, where $\F$ stands for $\R$ or $\C$. A linear
operator $N:\g \rightarrow \g$ is called Nijenhuis if
\begin{eqnarray}
\sabz\label{e10.10} [Nx,Ny]-N([Nx,y]+[x,Ny])+N^2[x,y]=0
\end{eqnarray}
for all $x,y\in\g$.
\end{defi}
The word "algebraic" in the titles of this section and the paper is used to
distinguish the algebraic situation from the geometric one,
where Nijenhuis operators are the endomorphisms of the
tangent bundle to a manifold. Since these last will not be used in the paper we shall omit the term algebraic.

Given a Nijenhuis operator, it can be showed (see \cite{mks})
that the operation 
\begin{equation}\sabz\label{e10.10.2}
[,]_N:\g\times\g \to\g,
[x,y]_N:=[Nx,y]+[x,Ny]-N[x,y],\ x, y\in\g, 
\end{equation}
is again a Lie
algebra bracket and moreover so is any linear combination
$[x,y]^\la:=[x,y]_N-\la[x,y],\la\in\F\cup\infty$. This facts
also follow from the proposition below, which will be used
for calculating the Casimir functions of  the Poisson pencil
corresponding to the pencil of Lie algebras $(\g,[,]^\la)$ and for other purposes.
\abz\label{10.20}
\begin{propo}
Let $N$ be a Nijenhuis operator. Then the following equality
holds for any $x,y\in \g$ and any $\la\not\in\Sp N$, where
$\Sp N$ stands for the spectrum of $N$:
\begin{equation}\sabz\label{e10.20}
(N-\la\Id)^{-1}[(N-\la\Id)x,(N-\la\Id)y]=[x,y]_N-\la[x,y]=[x,y]^\la.
\end{equation}
\end{propo}
\pf is straightforward:
\begin{eqnarray*}
(N-\la\Id)^{-1}[(N-\la\Id)x,(N-\la\Id)y]=(N-\la\Id)^{-1}([Nx,Ny]-
\la([Nx,y]+[x,Ny])+\la^2[x,y])&=&\\
(N-\la\Id)^{-1}(N([Nx,y]+[x,Ny])-N^2[x,y]-
\la([Nx,y]+[x,Ny])+\la^2[x,y])&=&\\
(N-\la\Id)^{-1}((N-\la\Id)([Nx,y]+[x,Ny]-N[x,y])-
\la(N[x,y]-\la[x,y]))&=&\\
(N-\la\Id)^{-1}((N-\la\Id)([x,y]_N-\la[x,y]))&=&[x,y]^\la.
\end{eqnarray*}
\qed

Now,  the LHS of the proved equality is a Lie bracket for
almost all $\la$, hence by continuity $[,]^\la$ is a Lie bracket for all $\la$.

The next lemma together with results of \cite{mks} allow to give a description of Nijenhuis operators, which is complete in the
diagonalizable case.
\abz\label{10.30}
\begin{lemm}
Let $N$ be a Nijenhuis operator. Then $N-\la\Id$ is a
Nijenhuis operator for any $\la\in\F\cup\infty$.
\end{lemm}

\pf Let's note that by definition an invertible operator $A$
is Nijenhuis if and only if $A^{-1}[Ax,Ay]=[x,y]_A$ for any
$x,y\in\g$. Now, by previous proposition for $A:=N-\la\Id$
and for almost all $\la$ we have
$$
A^{-1}[Ax,Ay]=[Nx,y]+[x,Ny]-N[x,y]-\la[x,y]=[x,y]_A.
$$
By continuity we conclude that $N-\la\Id$ is Nijenhuis for
any $\la$. \qed
\abz\label{10.35}
\begin{rema}\rm
One can also prove  that any "linear fractional" function $(s_1N+s_2\Id)(s_3N+s_4\Id)^{-1}, s_1,\ldots,s_4\in\F$, of a Nijenhuis operator is again a Nijenhuis operator.
\end{rema}
\abz\label{10.40}
\begin{propo}
Let $\g$ be a Lie algebra over $\C$, let $N$ be a Nijenhuis
operator, and let $\g=\g_1\oplus\cdots\oplus\g_n$ be its
decomposition to root spaces. Then this decomposition has the
following property: the subspace of the form
$\g_{i_1}\oplus\cdots\oplus\g_{i_k}, i_1<\cdots<i_k$, is a
Lie subalgebra for any $k<n$ (equivalently, $g_i+g_j$ is a
subalgebra for any $i,j$).

Conversely, any direct decomposition of $\g$ to subspaces
with the property above determines a diagonalizable Nijenhuis
operator uniquely up to a choice of eigenvalues
$\la_1,\ldots,\la_n$ corresponding to the subspaces
$\g_1,\ldots,\g_n$. In particular, any decomposition
$\g=\g_1\oplus\g_2$ to two subalgebras determines a Nijenhuis
operator.
\end{propo}

\pf Let $N$ be Nijenhuis and let $\Sp
N=\{\la_1,\ldots,\la_n\}$. By \cite[Subsection 2.1]{mks} and
by Lemma \ref{10.20} the subspaces $\g_i=\ker A_i^{r_i},
\check{\g}_i:=\oplus_{j\not=i}\g_j=\im A_i^{r_i}$ are Lie
subalgebras, where we put $A_i:=N-\la_i\Id$ and $r_i$ is the
Riesz index of $A_i$, i.e. the smallest integer with the
property that $\im A_i^{r_i}=\im A_i^{r_i+1}=\cdots$, while
$\im A_i^{r_i-1}\not=\im A_i^{r_i}$. Obviously, the
restriction of a Nijenhuis operator to a subalgebra is again
a Nijenhuis operator. So we can pass to $N|_{\check{\g}_i}$
and repeat the  considerations above. By induction we get the
desired property.

Now, let the decomposition $\g=\g_1\oplus\cdots\oplus\g_n$ be
such that $\g_{i}+\g_{j}$ is a
subalgebra for any $i,j$. Define $N$ by $N|_{g_i}:=\la_i\Id_{g_i},i=1,\ldots,n$. By the
bilinearity it is enough to prove equality (\ref{e10.10}) for
$x\in\g_i, y\in\g_j,1\le i,j\le n$:
\begin{eqnarray*}
[Nx,Ny]-N([Nx,y]+[x,Ny])+N^2[x,y] & =\\
\la_i\la_j([x,y]_i+[x,y]_j)-N(\la_i ([x,y]_i+[x,y]_j)+
\la_j([x,y]_i+[x,y]_j))+N(\la_i[x,y]_i+\la_j[x,y]_j) & =\\
\la_i\la_j([x,y]_i+[x,y]_j)-(\la_i^2[x,y]_i
+\la_i\la_j[x,y]_j+\la_i\la_j[x,y]_i
+\la_j^2[x,y]_j)+(\la_i^2[x,y]_i+\la_j^2[x,y]_j) & = & 0
\end{eqnarray*}
(here we denote by $[x,y]_i$ the i-th component of the
element $[x,y]$ with respect to the decomposition above).
\qed 

Some examples of Nijenhuis operators can be found in
\cite{mks}, another can be built using the second part of the
proposition above.

The fundamental example for this paper is as follows.
\abz\label{10.50}
\begin{exa}\rm
 Let $\g$ be an associative algebra and the Lie bracket is
 the commutator: $[x,y]:=xy-yx$. Then the operator $L_a$ of left
 (associative) multiplication by an element $a\in \g$ is a Nijenhuis
 operator: $[L_ax,L_ay]-L_a([L_ax,y]+[x,L_ay])+
 L_a^2[x,y]=axay-ayax-a(axy-yax+xay-ayx)
+a^2(xy-yx)=0$. In particular, if $\g=\gl_n$ is the algebra of
$n\times n$-matrices we get important examples of: a)
nilpotent Nijenhuis operator if $a\in \gl_n$ is nilpotent; b)
diagonalizable Nijenhuis operator if $a\in \gl_n$ is
diagonalizable. If $a=\diag(\la_1,\ldots,\la_n)\in \gl_n$ is
diagonal with $\la_i\not=\la_j$ while $i\not=j$, then the
corresponding eigenspaces $\g_i$ are equal to the Lie
subalgebras of matrices whose all rows except the i-th one
are the zero vectors.

It is easy to see that $[x,y]_{L_a}=xay-yax$ and the
corresponding pencil of Lie brackets is of the form
$[x,y]^\la=x(a-\la)y-y(a-\la)x$.
\end{exa}
Now  we want to remark that the last formula admits a generalization for a wider class of
Lie algebras. Below we shall give a
construction of pencils of Lie brackets on a class of
subalgebras of $\gl_n$, which also come from a Nijenhuis
operator but this operator in a sense is an "outer" one. 
\abz\label{10.60}
\begin{exa}\rm
Fix a matrix $I\in \gl_n$ and put $\g_I:=\{B\in \gl_n\mid
BI+IB^*=0\}, \h_I:=\{A\in \gl_n\mid AI-IA^*=0\}$, where $*$
denotes some involution on $\gl_n$ such that $(AB)^*=B^*A^*$
for any  $A,B\in \gl_n$. Then it is easy to see that $\g_I$ is
a Lie subalgebra in $\gl_n$ and that so is $L_A\g_I$ for any
$A\in \h_I$:
$[B,C]I=BCI-CBI=-BIC^*+CIB^*=IB^*C^*-IC^*B^*=-I([B,C])^*,
[AB,AC]=A(BAC-CAB),
(BAC-CAB)I=-BAIC^*+CAIB^*=-BIA^*C^*+CIA^*B^*=IB^*A^*C^*-IC^*A^*B^*
=-I(BAC-CAB)^*,B,C\in\g_I$. This shows that the formula
$[B,C]_A:=BAC-CAB$ defines a new Lie bracket on $\g_I$. Since
for any $\la\in\F$ we have $A-\la I_n\in \h_I$, where $I_n$
is the unity matrix, the brackets $[,]_A$ and $[,]$ generate
the pencil of Lie brackets $[,]^\la:=[,]_{(A-\la
I_n)}=[,]_A-\la[,]$. In general, this pencil is not generated
by an "inner" Nijenhuis operator because in general
$L_A\g_I\not=\g_I$. However, the formula of Proposition
\ref{10.20} is still valid: $[B,C]^\la=(A-\la
I_n)^{-1}[(A-\la I_n)B,(A-\la I_n)C], B,C\in\g_I$.
\end{exa}
\abz\label{10.70}
\begin{rema}\rm One can generalize this construction to the following one. Let $N$ be a Nijenhuis operator on $\g$ and let $\k \subset \g$ be a Lie subalgebra such that $(N-\la\Id)\k$ again is a Lie subalgebra for any $\la$. Then $(N-\la\Id)^{-1}[(N-\la\Id)\cdot,(N-\la\Id)\cdot]$ is a Lie bracket on $\k$ which is equal to $[\cdot,\cdot]_N|_\k-\la[\cdot,\cdot]|_\k=:[\cdot,\cdot]^\la$ by Proposition \ref{10.20}. In particular, $[\cdot,\cdot]^\la$ is a correctly defined pencil of Lie brackets on $\k$. 
\end{rema}
\section{Preliminaries on Poisson pencils and
formulation of main results}

All definitions below admit the real ($C^\infty$) and the complex
(holomorphic) versions. However, for the purposes of this
paper we shall need only the last one. So all
objects in the next two sections are complex analytic, $M$ stands for a
connected manifold. We refer
the reader to the book \cite{wcs} for the preliminaries on
Poisson structures.
\abz\label{20.10}
\begin{defi}\rm
A pair $(\t_1,\t_2)$ of linearly independent bivector fields
(bivectors for short) on a manifold $M$ is called Poisson if
$\t^s:=s_1\t_1+s_2\t_2$ is a Poisson bivector for any
$s=(s_1,s_2)\in \C^2$; the whole family of Poisson bivectors
$\{\t^s\}_{s\in\C^2}$ is called a Poisson pencil or a
bi-Poisson structure (or bihamiltonian structure).
\end{defi}
A bi-Poisson structure $\{\t^s\}$ (we shall often skip the
parameter space in the notations) can be viewed as a
two-dimensional vector space of Poisson bivectors, the
Poisson pair $(\t_1,\t_2)$ as a basis in this space. Of
course, the basis can be changed.
\abz\label{20.15}
\begin{exa}\rm
Let $\g$ be a Lie algebra over $\C$ with a Nijenhuis operator
$N$. Denote by $\t_1,\t_2$ the canonical linear Poisson
bivectors (the so-called Lie--Poisson bivectors) on the dual
space $\D$ related to the Lie brackets $[,]$ and $[,]_N$,
respectively. Then, since these brackets generate a pencil of
Lie brackets (see Section \ref{s10}), the pair $\t_1,\t_2$ is
Poisson. The corresponding Poisson pencil will be denoted by
$\{\t_N^s\}$. For a real case we complexify all objects and then build the
holomorphic Poisson pencil on $(\g^\C)^*$ as above.
\end{exa}

The following definition is due to I.Zakharevich \cite{z1}
\abz\label{20.20}
\begin{defi}\rm
Let $\{\t^s\}$ be a Poisson pencil on $M$. It is called
Kronecker at a point $x\in M$ if $\rk_\C\t_x^s$ is constant
with respect to $s\in\C^2\setminus\{0\}$. We say that
$\{\t^s\}$ is micro-Kronecker (Kronecker for short) if it is
Kronecker at any point of some open dense set in $M$.
\end{defi}

The next proposition shows the importance of Kronecker
Poisson pencils, which serve as a convenient formalism
allowing to construct and investigate completely integrable
systems. For the proof see \cite{bols}, \cite{p1}.
\abz\label{20.30}
\begin{propo}
Let $\{\t^s\}$ be a Kronecker Poisson pencil on $M$. Assume
that an open set $U\subset M$ is such that the set
$Z^{\t^s}(U)$ of Casimir functions for $\t^s$ over $U$ is
complete (i. e. the common level sets of functions from
$Z^{\t^s}(U)$ coincide with the regular part of the
symplectic foliation of $\t^s$ on $U$) for any $s\not=0$.
Then the set
\[
Z^{\{\t^s\}}(U):=\sum_{s\not=0}Z^{\t^s}(U)
\]
is a complete involutive set of functions for any $\t^s\not
=0$, that is, the common level sets of functions from
$Z^{\{\t^s\}}(U)$ form a lagrangian foliation in any regular
symplectic leaf of $\t^s$ on $U$.
 (Here  we understand the sum as the linear span of an
infinite family of linear subspaces of functions, which in
fact is generated by some finite subfamily.)
\end{propo}
We shall call the functions from  $Z^{\{\t^s\}}$ the first
integrals of the Poisson pencil $\{\t^s\}$.

Now we are ready to formulate the main result of this paper
which gives necessary and sufficient conditions of the
kroneckerity of the Poisson pencil $\{\t_N^s\}$ built on a
Nijenhuis operator $N$ (see Example \ref{20.15}). 

Given a homogeneous space $G/H$, where $G\supset H$ are Lie groups, the Lie group $H$ is naturally acting on it and since the point $eH$, where $e\in G$ is the neutral element, is stabilized by this action one can extend it to the linear action in the cotangent space $T^*_{eH}(G/H)$. This representation of $H$ is denoted by $\rho$ and is called the coisotropy representation. Identifying $T_{eH}(G/H)$ with $\g/\h$ we obtain the following formula:
$\rho:H\to \Aut((\g/\h)^*), h\stackrel{\rho}{\mapsto}
(g+\h\mapsto\Ad_h(g)+\h)^*,g\in\g,h\in H$. For any element $a\in(\g/\h)^*$ we introduce two numbers: $\ind a$, which is the index of the Lie algebra $\h^a$ of the stabilizer $H^a$ of $a$ with respect to the coisotropy action, and $\codim a$, which is the  codimension of the orbit $H\cdot a$ of $a$ with respect to the coisotropy action (recall that the index of a Lie algebra is by definition the codimesion of the generic coadjoint orbit). 

\abz\label{main1}
\begin{theo}
Let $\g$ be a Lie algebra over $\C$ and let $N$ be a 
diagonalizable Nijenhuis operator with the spectrum $\Sp
N=\{\la_1,\ldots,\la_n\}$. Assume the following condition is satisfied:
\begin{equation}\sabz \label{e.main}
{the\ complement\ to\ the\ set} \bigcup_{\la\in\C\setminus\Sp N}(N-\la\Id)^{-1}(\Sing\g^*)\ in\ \g^*\ contains\ an\ open\ dense\ set.
\end{equation}
 (Here we have denoted by $\Sing\g^*$ the union of all coadjoint orbits of nonmaximal dimension. In particular, if $\g$ is reductive, $\codim\Sing\g^*\ge 3$ and the assumption above is satisfied.)

Put $\check{\g}_i:=\im(N-\la_i\Id),i=1,\ldots,n$ ($\check{\g}_i$ are Lie
algebras by Lemma \ref{10.30} and by equality (\ref{e10.10}), see also Proposition \ref{10.40}).

Then the corresponding Poisson
pencil $\{\t_N^s\}$ is Kronecker if and only if  one can find elements $c_1,\ldots,c_n, c_i\in (\g/\check{\g}_i)^*, i=1,\ldots,n$, such that for any $i,1\le i\le n$,
\begin{equation}\sabz\label{e.main1}
\ind c_i+\codim c_i=\ind \g,
\end{equation}
where $\ind c_i,\codim c_i$ are the corresponding numbers related to the coisotropy representation $\rho_i:\check{G}_i\to (\g/\check{\g}_i)^*$, $\check{G}_i$ being the Lie subgroup of $G$ with the Lie algebra $\check{\g}_i$.
\end{theo}

The proof of this result is postponed to
Section \ref{proofs}. Taking $c_i=0,i=1,\ldots,n$, we get the following corollary.
\abz\label{coro}
\begin{coro}
Under the assumptions of Theorem \ref{main1}, if for any $i, 1\le i\le n$,
$$
\ind\check{\g}_i+\codim\check{\g}_i=\ind\g,
$$
then the Poisson pencil $\{\t_N^s\}$ is Kronecker.
\end{coro}

The last part of this section is devoted to some definitions
which are based on the reference \cite{mks88} and which will be used later. 

Let $\g=\g_1\oplus\g_2$ be a Lie algebra which is a direct sum of its subalgebras $\g_1,\g_2$. Then the Lie bracket on $\g$ can be decomposed as follows: 
\begin{equation}\sabz\label{e20.40.0}
[x,y]=[x_1,y_1]_1+([x_1,y_2]_1+[x_2,y_1]_1)+([x_1,y_2]_2+[x_2,y_1]_2) +[x_2,y_2]_2,
\end{equation}
where the indices refer to the corresponding projections onto $\g_1$ or $\g_2$. It turns out that the maps $A_1\colon x_1\mapsto [x_1,\cdot]_2\colon \g_1\to \End(\g_2)$ and $A_2\colon [x_2,\cdot]_1\colon \g_2\to \End(\g_1)$ are the Lie algebra homomorphisms, where $\End(\g_i)$ is a Lie algebra of the endomorphisms of the vector space $\g_i$. The representations $A_1,A_2$ also satisfy some additional conditions making them cocycles, which will be inessential for us. The formula above rewritten in terms of the pairs of elements $(x_1,x_2),x_i\in\g_i$, is the following:
\begin{equation}\sabz\label{e20.40.1}
[(x_1,x_2),(y_1,y_2)]=([x_1,y_1]+ A_2(x_2)y_1-A_2(y_2)x_1,[x_2,y_2]+A_1(x_1)y_2-A_1(y_1)x_2).
\end{equation}
In particular, if one put here $A_2\equiv 0$, one gets
the formulas
\begin{equation}\sabz\label{e20.40.2}
[(x_1,x_2),(y_1,y_2)]=([x_1,y_1],A_1(x_1)y_2-A_1(y_1)x_2),
\end{equation}
or
\begin{equation}\sabz\label{e20.40.3}
[x,y]=[x_1,y_1]_1+([x_1,y_2]_2+[x_2,y_1]_2),
\end{equation}
in which one recognizes the  multiplication in the semidirect product $\g_1\times_{A_1}\g_2$, where $\g_2$ is regarded as a vector space (abelian Lie algebra).  
\abz\label{20.40}
\begin{defi} 
We shall refer to the direct sum of Lie algebras $\g_1\oplus\g_2$ with the bracket (\ref{e20.40.0}) (or (\ref{e20.40.1})) as to a twilled Lie algebra and shall denote it by $\g_1\bowtie\g_2$. The same vector space with the multiplication (\ref{e20.40.2}) (or (\ref{e20.40.3})) will be called the truncation of the twilled Lie algebra $\g_1\bowtie\g_2$ and denoted by $\g_1\sdp\g_2$. 
\end{defi}

\section{Proof of the main result}
\label{proofs}
\noindent{\sc Proof of Theorem \ref{main1}} 
The idea of the proof is as follows. As we have already mentioned in the Introduction formula (\ref{e10.20}) shows that almost all bivectors of the Poisson pencil $\{\t_N^s\}$ built from the algebraic Nijenhuis operator $N$ are isomorphic, in particular have the same corank equal to the index of the algebra $\g$. The exception are the bivectors $\t^{(-\la_i,1)}, i=1,\ldots,n$, corresponding to the eigenvalues $\la_i$ of $N$ and all that we need to control the corank of the whole pencil is to control the rank of these bivectors which will be called exceptional. To this end we shall show that the Lie brackets corresponding to the exceptional bivectors are semidirect products and we shall use the Ra\"{\i}s type formula for the index of such algebras.

Let us fix $i, 1\le i\le n$, and consider the Nijenhuis operator $M:=N-\la_i\Id$. Then the following formula describes the deformed bracket corresponding to $M$ (see \cite[Section 6]{cgm1}):
$$
[x,y]_M=M|_{E_1}^{-1}[Mx,My]_1+[Mx,y]_2+[x,My]_2,
$$
where we put for a moment $E_1:=\check{\g}_i=\im M,E_2:=\g_i=\ker M$ and the subscripts refer to the projections onto $E_1$ or $E_2$. We claim that the bracket $[,]_M$ is the truncated bracket (see Definition \ref{20.40}) corresponding to the twilled Lie algebra structure on $\g=E_1\times E_2$ given by the bracket $[,]_L$ defined below. 

Define a new Nijenhuis operator $L$ on $\g$ by the formula $L=M\circ P_1+P_2$, where $P_i$ stands for the projector onto $E_i,i=1,2$. In other words $L$ acts as $M$ on $E_1$ and identically on $E_2$. Now let $[,]_L$ stand for the deformed bracket corresponding to $L$ by formula (\ref{e10.10.2}). Then  by (\ref{e10.20}) we have
\begin{eqnarray*}
[x,y]_L=L^{-1}[Lx,Ly]=L|_{E_1}^{-1}P_1[Lx,Ly]+P_2[Lx,Ly]\\ =
 M|_{E_1}^{-1}([Mx_1,My_1]_1+[Mx_1,y_2]_1+[x_2,My_1]_1) + 
[Mx_1,y_2]_2+[x_2,My_1]_2+[x_2,y_2]_2.
\end{eqnarray*}
Now, the truncated bracket equals $M|_{E_1}^{-1}[Mx_1,My_1]_1+
[Mx_1,y_2]_2+[x_2,My_1]_2$ which coincides with $[x,y]_M$ and the claim is proved.

Note that the truncated brackets corresponding to the twilled Lie algebras $(\check{\g}_i\times\g_i,[,])$ and $(\check{\g}_i\times\g_i,[,]_L)$ are isomorphic. Indeed, the isomorphism is given by the operator $L:(\check{\g}_i\times\g_i,[,]_L)\to (\check{\g}_i\times\g_i,[,])$. This isomorphism is compatible with the truncations, i. e. the corresponding truncated algebras also are isomorphic. In particular, in the considerations below concerning the codimensions of the coadjoint orbits we can regard simply the semidirect products which are truncations of the twilled Lie algebras $(\check{\g}_i\times\g_i,[,]), i=1,\ldots,n$. 

Now we can use the standard facts about semidirect products, which we recall below. Given a semidirect product $\g\times_\rho V$ of a Lie algebra $\g$ with a vector space $V$ by means of a representation $\rho$, one can show 
(see \cite{rst}, for example) that, 
1) any covector $a\in V^*$ is contained in a  set $V_a \subset (\g\times_\rho V)^*$ which is a Poisson submanifold in $(\g\times_\rho V)^*$ isomorphic to $T^*G/G^a$ (here $G,G^a$ are the Lie groups corresponding to the Lie algebras $\g,\g^a$, where $\g^a$ is the stabilizer of $a$); 2) the coadjoint orbits contained in $V_a$ are isomorphic to the symplectic leaves of $T^*G/G^a$, in particular, the generic (in $V_a$) orbits has codimension in $V_a$ equal to $\ind \g^a$; 3) the submanifold $V_a$ is of the form $\g^*\times O_a$, where $O_a \subset V^*$ is the orbit of $a$ in $V^*$, in particular, $\codim_{(\g\times_\rho V)^*}V_a=\codim_{V^*}O_a$. Summarizing all this, we can say that, given an element $a\in V^*$,  one can associate  with it a coadjoint orbit $S_a$ of  $(\g\times_\rho V)^*$ such that $\codim_{(\g\times_\rho V)^*}S_a=\codim_{V^*}O_a+\ind\g^a$. Taking a generic $a$ we get the so called Ra\"{\i}s formula \cite{rais}: $\ind (\g\times_\rho V)^*=\codim_{V^*}O_a+\ind\g^a$.

Now let us complete the proof of Theorem \ref{main1}. We claim that 
the kroneckerity of the Poisson pencil $\{\t^s_N\}$ is equivalent to existing for  any $i,1\le i\le n$, of a symplectic leaf $S_i$ of the bivector $\t_N^{(-\la_i,1)}$ such that $\codim S_i=\ind\g$.
Indeed, by formula (\ref{e10.20}) we have $(N-\la\Id)^{-1}_*\t_1=\t_N^{(-\la,1)}$ for any $\la\not\in\Sp N$. Since $\Sing\g^*$ is the union of symplectic leaves of $\t_1$ of nonmaximal dimension, the assumption (\ref{e.main}) implies that on the open dense set mentioned in it the corank of the bivectors $\t_N^{(-\la,1)},\la\not\in\Sp N$, is equal to $\ind\g$. But these bivectors up to rescaling exhaust all nonexceptional bivectors. Now it is clear that the kroneckerity is equivalent to the condition
$\crk\t_N^{(-\la_i,1)}=\ind\g, i=1,\ldots,n$. On the other hand, in general $\ind\g\le\crk\t_N^{(-\la_i,1)}\le\codim S_i$ and we have proved the claim. 

Now let us pick out a point $c_i\in(\g/\check{\g}_i)^*\simeq(\g_i)^*$ for any $i=1,\ldots,n$. By the considerations above this is equivalent to picking out a symplectic leaf $S_i$ of the exceptional bivector $\t^{(-\la_i,1)}$ such that $\codim_{\g^*}S_i=\codim_{\g_i^*}O_{c_i}+\ind\check{\g}_i^{c_i}$, where $O_{c_i}$ is the $\check{G}_i$-orbit of $c_i$ in $\g_i^*$ and $\check{\g}_i^{c_i}$ is its stabilizer. \qed

\section{First example: Kronecker Poisson pencils related to Borel subalgebras}
\label{s.borel}

Let $\g$ be a semisimple real or complex Lie algebra of rank $r$. If $\g$ is real, assume that it is split. Consider a Cartan sublagebra $\h$ and the corresponding direct decompositions 
$\g =\n_-\oplus\b=\n_-\oplus\h\oplus\n_+$, where $\b,\n_-,\n_+$ are the Borel and the maximal nilpotent subalgebras respectively.
Define a linear operator $N: \g\to\g$ by $N|_{\n_-}=\la_1\Id, N|_{\b}=\la_2\Id$. Then by Proposition \ref{10.40} $N$ is a Nijenhuis operator.
\abz\label{borel}
\begin{theo}
The corresponding Poisson pencil $\{\t^s_N\}$ is Kronecker.
\end{theo}

\pf We need to check that $(\g,N)$ satisfies the criterion of the kroneckerity, Theorem ref{main1}. Note that since $\g$ is semisimple, $\codim\Sing\g\ge 3$ and condition (\ref{e.main}) is satisfied and also $\ind\g=\rk\g$. 

Now we shall consider the corresponding  coisotropy representations. Denote by $N_-,B$ the corresponding Lie groups. Using the Killing form we obtain the following natural identifications: $(\g/\b)^*\simeq\n_+, (\g/\n_-)^*\simeq\b_-:=\n_-\oplus\h$, which are $B$- and $N_-$-equivariant respectively. We need to find two elements $c_1\in\n_+$ and $c_2\in\b_-$ satisfying condition (\ref{e.main1}).

First take $c_1:=e$, a principal nilpotent element \cite{kost,bourb}. Its stabilizer $\g^e$ with respect to the adjoint action of $\g$
is an abelian subalgebra of dimension $r$. Moreover, $\g^e \subset \n_+$, hence $\g^e=\b^e$ and $\ind c_1=r$. The dimension of the adjoint $B$-orbit of $c_1$ equals $\dim\n_+$, i. e. $\codim c_1=0$ and (\ref{e.main1}) is satisfied for $c_1$.

Now let $c_2:=f\in\h$  be a regular semisimple element. Then  $\g^f=\h$ and $\n_-^f=\g^f\cap\n_-=0$. Thus $\ind c_2=0$ and the dimension of the $N_-$-orbit of $c_2$ in $\b_-$ is $\dim\n_-$, i. e. $\codim c_2=\dim\b_--\dim\n_-=r$. \qed
\abz\label{borel1}
\begin{rema}\rm
Let us exhibit what functions in involution we get on $\g^*$ with the help of $N$ (for general pair $(\g,N)$). Recall that they are generated by the Casimir functions of all bivectors of the pencil $\{\t^s_N\}$ (see Proposition \ref{20.30}).By formula (\ref{e10.20}) the Casimirs of $\t^{s}_N, s=(s_1,s_2)$, are functionally generated by $C_j((N-\la\Id)^{-1}x),j=1,\ldots,r$, where $\la=-s_1/s_2$ and $C_i,\ldots,C_r$ are the independent Casimirs of $\t_N^{(0,1)}$, i. e. the invariants of the coadjoint action of $\g$.
In fact it is enough to choose a finite number of bivectors $\t^{s_i}_N, s_i=(s_i^1,s_i^2),i=1,\ldots,p$, where $p$ is sufficiently large, and one can take $-s_i^1/s_i^2$ to be not equal to the eigenvalues of $N$. Thus our family of functions in involution is functionally generated by $C_j((s_i^2N+s_i^1\Id)^{-1}x),j=1,\ldots,r, i=1,\ldots,p$.

If $C_1,\ldots,C_r$ are polynomials as in the case of semisimple $\g$, another way to obtain this family of functions is to consider the coefficients of the expansion of $C_j((N-\la\Id)^{-1}x)$ in the powers of $\la$ (or $1/\la$).
\end{rema}

Taking $\la_1=1,\la_2=-1$ for the Nijenhuis operator $N$ built above we obtain the following simple formula for the resolvent: $(N-\la\Id)^{-1}=(1-\la^2)^{-1}(N+\la\Id)$. In particular, for $\g=\sl _n$ we can take the coefficients of the expansion in $\la$ of the following functions $\Tr((N+\la\Id)x)^k, k=2,\ldots,n$ (we have identified $\g$ and $\g^*$ by means of the Killing form). These functions form a complete involutive set on any adjoint orbit of maximal dimension. It is easy to see that the following quadratic hamiltonians are in this family: $\sum_{i=1}^nx_{ii}^2, \sum_{i<j}x_{ij}x_{ji}$.

\section{Second example: $n$-dimensional free rigid body and the method of argument translation}
\label{s.manakov}
\abz\label{main2}
\begin{theo}
Let $\g=\gl_n$ be the Lie algebra of $n\times n$-matrices and
let $N=L_A$ be a Nijenhuis operator of left multiplication by
the diagonal matrix $A=\diag(\la_1,
\ldots,\la_n)\in\g$
with $\la_i\not=\la_j,i\not=j$.

Then the corresponding Poisson pencil $\{\t_N^s\}$ is
Kronecker.
 \end{theo}

\pf Obviously, the subalgebra $\check{\g}_i$ for the Nijenhuis operator $L_A$ equals the set of matrices with the zero i-th row, hence $\codim\check{\g}_i=n$.  
The proof follows from Corollary \ref{coro} and from the following lemma. \qed
\abz\label{proofs.20}
\begin{lemm}
Let $\check{\g}_i \subset \g=\gl_n$ be the Lie subalgebra of matrices with the zero i-th row. Then $\ind \check{\g}_i=0$ (i.e. $\check{\g}_i$ is Frobenius).
\end{lemm}

\pf To prove the lemma we shall use the fact that $\t_i=p_*\t$, where $\t_i,\t$ are the Lie--Poisson bivectors on $\check{\g}_i^*$ and $\g^*$ correspondingly, and  $p:\g^*\to\check{\g}_i^*=\g^*/\check{\g}_i^\bot$ is the natural projection. Using the nondegenerate invariant form $\Tr(XY)$ we identify $\g^*$ with $\g$ and the annihilator $\check{\g}_i^\bot$ with the space $V_i$ of matrices with all but i-th zero columns. 

The following assertion implies the lemma: the intersection of the  adjoint orbit $O$ of a generic (equivalently, of some) element $X\in\g$ with the affine subspace $X+V_i$ is trivial. Indeed, the map $p$ is Poisson as a map between the Lie--Poisson structures and,  since the dimensions of the orbit $O$ and the subspace $X+V_i$ are complementary, the triviality of the intersection implies that codimension of the orbit $p(O)$ is zero.

Now we shall  prove the above mentioned assertion.  Let $M=M(m_1,\ldots,m_n)\in V_i$ be a matrix whose i-th column constitute the elements $m_1,\ldots,m_n$. Assume  $X$ is generic semisimple and  $X+M\in O$. Then $\Tr(X+M)^k=\Tr X^k,k=1,\ldots,n$. The first of this equalities implies $m_i=0$ which yields $M^k=0$ for any $k>1$. Thus we get $\Tr(X+M)^k=k\Tr(X^{k-1}M)+\Tr X^k$, hence $\Tr(X^{k-1}M)=0, k=2,\ldots,n$. In other words we have the system of linear equations
$$
x^{(k)}_{1i}m_1+\cdots+x^{(k)}_{i-1,i}m_{i-1}+x^{(k)}_{i+1,i}m_{i+1}+
\cdots+x^{(k)}_{ni}m_n=0,\  k=1,\ldots,n-1,
$$
where $x^{(k)}_{jl}$ are the coefficients of the matrix $X^k$. Now it remains to show that $X$ can be chosen so that the vectors $w_k(X):=(x^{(k)}_{1i},\ldots,x^{(k)}_{i-1,i},x^{(k)}_{i+1,i},\ldots,x^{(k)}_{ni}),
k=1,\ldots,n-1$, are linearly independent. 

Let's take
$X:=\diag(a_1,\ldots,a_{i-1},0,a_{i+1},\ldots,a_n)+L_i$, where $a_j$ are different not equal to 1 and $L_i$ is the matrix whose all entries are zero except that from the i-th row which are equal to 1. Obviously, the eigenvalues of $X$ are different, i.e. $X$ is generic semisimple. Also it is easy to calculate that the i-th row of $X^k$ equals 
$(\sum_{j=0}^{k-1}a_1^j,\ldots,\sum_{j=0}^{k-1}a_{i-1}^j,1,
\sum_{j=0}^{k-1}a_{i+1}^j,\ldots,\sum_{j=0}^{k-1}a_n^j)$. So the vectors $w_1(X),\ldots,w_{n-1}(X)$ are linearly independent since their linear span coincides with that of the columns of the corresponding Vandermonde matrix. \qed
 
Let us look at the functions in involution obtained from this example. By Remark \ref{borel1} they can be functionally generated by the coefficients of the expansion of the functions $\Tr((N-\la\Id)^{-1}x)^k, k=1,\ldots,n$, in $\la$ (we have identified $\g$ and $\g^*$ by means of the "trace" form). Let us rewrite the resolvent $(N-\la\Id)^{-1}$ in the form $-\la^{-1}(\Id-\frac{1}{\la}N)^{-1}=
-\la^{-1}\sum_{j=0}^\infty\frac{1}{\la^j}N^j$. Write $f_{kl}$ for 
the coefficient of $1/\la^l$ in $\Tr(\sum_{j=0}^\infty\frac{1}{\la^j}N^jx)^k$.
Recall also that the so-called Manakov integrals $h_{kl}$ are the coefficients of $\la^{l}$ in $h_k^\la(x):=\frac{1}{k}\Tr(x+\la A)^k$.
\abz\label{manakov}
\begin{propo}
The functions $h_{kl},k=1,\ldots,n,\ l=0,\ldots,k-1$, and $f_{kl},k=1,\ldots,n,\ l=0,\ldots,k-1$, generate the same families of functions in involution. 
\end{propo}

\pf We shall use the recursion relations satisfied by both the families.
\abz\label{manakov.lem}
\begin{lemm}
Let $\t_1,\t_2$ be the Lie-Poisson structures corresponding to the Lie brackets $[,],[,]_N$ on $\g$ respectively.
Then 1) $\t_1(h_{k+1,l+1})=\t_2(h_{kl})$; 2) $\t_1(f_{k,l+1})=\t_2(f_{kl})$.
\end{lemm}

\pf To prove the first relation we adapt the proof of the analogous relation for $\g=\so_n$ in \cite{morosip}. We notice that $d_xh_k^\la(x)=(x+\la A)^{k-1}$ (here the matrix in the RHS is a functional via the "trace" form) and $\t_1|_xx'=[x,x'],\t_2|_xx'=xx'A-Ax'x$ (here  $\t_i$ is a map $T^*\g^*\to T\g^*\simeq\g^*\stackrel{\Tr}{\simeq}\g$, $x'\in T^*_x\g^*\simeq\g$). Now it is straightforward to show that $\t_1h_{k+1}^\la=\la\t_2dh_k^\la$, which proves 1).

To show the second relation we use the fact that $\Tr((N-\la\Id)^{-1}x)^k$ is a Casimir for $\t_2-\la_1\t_2$, i. e. $(\t_2-\la_1\t_2)(\sum_j \frac{1}{\la_j}f_{kj})=0$. Comparing the coefficients we get 2). \qed

Now we are ready to prove Proposition \ref{manakov}. 
We have the following two hierarchies of functions:
$$
\begin{array}{cccc}
h_{10} &        &        &       \\
h_{20} & h_{21} &        &      \\
h_{30} & h_{31} & h_{32} &      \\
\vdots & \vdots & \vdots &\ddots
\end{array}\ \ \ \ 
\begin{array}{cccc}
   f_{10}&      &      &      \\
   f_{20}&f_{11}&      &      \\
   f_{30}&f_{21}&f_{12}&      \\
   \vdots&\vdots&\vdots&\ddots
\end{array}   
$$
We shall use the induction on the number $l$ of the column of the second table. It is easy to see that $h_{k0}=\frac{1}{k}\Tr x^k=\frac{1}{k}f_{k0}$. Now fix $l$ and  suppose  the function $f_{ml},1\le m\le n$, can be expressed as a function of $h_{kj}, j=0,\ldots,l,k=1,\ldots,n$. We should prove that $f_{m,l+1},1\le m\le n$, also can be expressed as a function of $h_{kj},j=0,\ldots, l+1,k=1,\ldots,n$. Indeed, since by the lemma above $\t_1(f_{m,l+1})=\t_2(f_{ml})$, this vector field can be expressed as a linear combination of vector fields $\t_2(h_{kj}),j\le l$, which in turn are equal to $\t_1(h_{k+1,j+1})$. This shows that the function $f_{m,l+1}$ can be functionally expressed by the functions $h_{kj},j=0,\ldots,l+1,k=1,\ldots,n$ (the Casimirs of $\t_1$ are included in this family). \qed

\bibliography{Main}

\providecommand{\bysame}{\leavevmode\hbox to3em{\hrulefill}\thinspace}
\providecommand{\MR}{\relax\ifhmode\unskip\space\fi MR }
\providecommand{\MRhref}[2]{%
  \href{http://www.ams.org/mathscinet-getitem?mr=#1}{#2}
}
\providecommand{\href}[2]{#2}
\begin{thebibliography}{RSTS94}

\bibitem[Bol91]{bols}
Alexey Bolsinov, \emph{Compatible {P}oisson brackets on {L}ie algebras and
  completeness of families of functions in involution}, Izv.Akad. Nauk SSSR,
  Ser. Mat. \textbf{55} (1991), In Russian.

\bibitem[Bou82]{bourb}
Nicolas Bourbaki, \emph{Groups et alg{\` e}bres de {L}ie, {IX}}, Masson, 1982.

\bibitem[CGM01]{cgm1}
Jos\'{e} Carin{\~e}na, Janusz Grabowski, and Giuseppe Marmo,
  \emph{Contractions: Nijenhuis and {S}aletan tensors for general algebraic
  structures}, J. Phys. A \textbf{34} (2001), 3769--3789.

\bibitem[dSW99]{wcs}
Ana~Canas da~Silva and Alan Weinstein, \emph{Geometric models for
  noncommutative algebras}, Providence, 1999.

\bibitem[GZ89]{gz1}
Israel Gelfand and Ilya Zakharevich, \emph{Spectral theory for a pair of
  skew-symmetrical operators on $s^1$}, Funktsion. Analiz i ego Prilozh.
  \textbf{23} (1989), no.~2, 1--11, In Russian.

\bibitem[GZ91]{gz2}
\bysame, \emph{Webs, {V}eronese curves, and bihamiltonian systems}, J. Funkt.
  Anal. \textbf{99} (1991), 150--178.

\bibitem[GZ93]{gz3}
\bysame, \emph{On the local geometry of a bihamiltonian structure}, The Gelfand
  mathematical seminars 1990-1992, Birkhauser, 1993, pp.~51--112.

\bibitem[GZ00]{gz4}
\bysame, \emph{Webs, {L}enard schemes, and the local geometry of bihamiltonian
  {T}oda and {L}ax structures}, Selecta-Math. (N.S.) \textbf{6} (2000),
  131--183.

\bibitem[Kos59]{kost}
B.~Kostant, \emph{The principal thee-dimensional subgroup and the {B}etti
  numbers of a complex simple {L}ie group}, Amer.J.Math. \textbf{81} (1959),
  973--1032.

\bibitem[KSM88]{mks88}
Yvette Kosmann-Schwarzbach and Franco Magri, \emph{Poisson-{L}ie groups and
  complete integrability {I}. {D}rinfeld bigebras, dual extensions and their
  canonical representations}, Ann. Inst. Henri Poincar\'{e} \textbf{49} (1988),
  433--460.

\bibitem[KSM90]{mks}
\bysame, \emph{Poisson-{N}ijenhuis structures}, Ann. Inst. Henri Poincar\'{e}
  \textbf{53} (1990), 35--81.

\bibitem[Man76]{man}
S.~V. Manakov, \emph{A remark on the integration of the {E}uler equation of the
  dynamics of an n-dimensional rigid body}, Funktsional. Anal. i Prilozhen.
  \textbf{10} (1976), no.~4, 93--94.

\bibitem[MF78]{mf}
A.~S. Mishchenko and A.~T. Fomenko, \emph{Euler equations on finite dimensional
  {L}ie groups}, Math. USSR, Izvestija \textbf{12} (1978), 371--389,
  Translation from the Russian.

\bibitem[MP96]{morosip}
Carlo Morosi and Livio Pizzocchero, \emph{On the {E}uler equation:
  Bi-{H}amiltonian structure and integrals in involution}, Lett. Math. Phys.
  \textbf{37} (1996), 117--135.

\bibitem[Pan00]{p1}
Andriy Panasyuk, \emph{Veronese webs for bihamiltonian structures of higher
  corank}, Banach Center Publications \textbf{51} (2000), 251--261.

\bibitem[Ra{\"\i}78]{rais}
Mustapha Ra{\"\i}s, \emph{L'indice des produits semi-directs
  ${E}\times_\rho\g$}, C. R. Acad. Sc. Paris, Ser. A \textbf{287} (1978),
  195--197.

\bibitem[RSTS94]{rst}
A.~G. Reyman and M.~A. Semenov-Tian-Shansky, \emph{Group-theoretical methods in
  the theory of integrable systems}, Encyclopaedia of Math. Sciences (Dynamical
  Systems VII), vol.~16, Springer, 1994, pp.~116--225.

\bibitem[STS94]{sts}
M.~A. Semenov-Tian-Shansky, \emph{Lectures on $r$-matrices, {P}oisson--{L}ie
  groups and integrable systems}, Lectures on {I}ntegrable systems, In memory
  of J.-L. Verdier, Proc. of the CIMPA School on {I}ntegrable systems, Nice
  1991 (Y.~Kosmann-Schwarzbach O.~Babelon, P.~Cartier, ed.), World
  {S}cientific, 1994, pp.~269--317.

\bibitem[Zak01]{z1}
Ilya Zakharevich, \emph{Kronecker webs, bihamiltonian structures, and the
  method of argument translation}, Transformation Groups \textbf{6} (2001),
  267--300.

\end{thebibliography}
\end{document}